\documentclass{acm_proc_article-sp}   
\usepackage{url,graphicx,array,colortbl}    
\usepackage{rotating}

\newtheorem{theorem}{Theorem}[section]
\newtheorem{example}[theorem]{Example}

\newfont{\saffaddrit}{phvro at 8pt} 
\newfont{\seaddfnt}{phvr at 9pt}
\def\semail#1{{{\seaddfnt{\vskip 4pt#1}}}} 

\title{Parallel computation of the rank of large sparse matrices
  from algebraic K-theory}
\numberofauthors{2}
\author{
  \alignauthor Jean-Guillaume Dumas\\[1mm]
  \medskip
  \saffaddrit{
    \mbox{Laboratoire Jean Kuntzmann,}
UMR CNRS 5224      \mbox{Universit\'e Joseph Fourier,}
    \mbox{B.P. 53 X, 38041 Grenoble, France.}\\
  }
  \medskip
  \semail{Jean-Guillaume.Dumas@imag.fr}
  \medskip
\and 
\alignauthor Philippe Elbaz-Vincent\\[1mm]
  \medskip
  \saffaddrit{
     Institut de Math\'ematiques et de Mod\'elisation de Montpellier,
UMR CNRS 5149\\
 \mbox{Universit\'e Montpellier II,}
 \mbox{CC051, Place E. Bataillon.}
 \mbox{34095 Montpellier cedex 5, FRANCE.}\\
  }
  \medskip
  \semail{pev@math.univ-montp2.fr}
  \medskip
\and
\alignauthor Pascal Giorgi\\[1mm]
  \medskip
  \saffaddrit{
     \mbox{Laboratoire LP2A,}
     \mbox{Universit\'e de Perpignan Via Domitia.}
     \mbox{52, avenue Paul Alduy}
     \mbox{66860 Perpignan France.}\\
}
  \medskip
  \semail{pascal.giorgi@univ-perp.fr}
  \medskip
\and 
\alignauthor Anna Urba\'nska\\[1mm]
  \medskip
  \saffaddrit{
    \mbox{Laboratoire Jean Kuntzmann,}
UMR CNRS 5224      \mbox{Universit\'e Joseph Fourier,}
    \mbox{B.P. 53 X, 38041 Grenoble, France.}\\
  }
  \medskip
  \semail{Anna.Urbanska@imag.fr}
  \medskip
}
\date{\today}

\newcommand{\Z}{\mathbb{Z}}
\newcommand{\Q}{\mathbb{Q}}
\newcommand{\R}{\mathbb{R}}
\newcommand{\N}{\mathbb{N}}

\renewcommand{\le}{\leqslant}

\definecolor{mongris}{rgb}{0.9, 0.9, .9}

\newcommand{\F}{{\sf F}} 
\def\sigmabase{$\sigma$-basis~}
\newcommand{\graycell}[1]{\multicolumn{1}{|>{\columncolor{mongris}}c|}{#1}}

\DeclareMathOperator{\rank}{\operatorname{rank}}

\begin{document}
\thispagestyle{empty}
\maketitle 
\thispagestyle{empty}
\pagestyle{empty}

\begin{abstract}
This paper deals with the computation of the rank and some integer Smith
forms of a series of sparse matrices arising in algebraic K-theory.
The number of non zero entries in the considered matrices ranges from
8 to 37 millions. 
The largest rank computation took more than 35 days on
50 processors. We report on the actual algorithms we used to build the
matrices, their link to the motivic cohomology and the linear algebra
and parallelizations required to perform such huge computations.
In particular, these results are part of the first computation of the cohomology of the linear group
$GL_7(\Z)$.
\end{abstract}

\section{Introduction}
\subsection{Motivation from K-theory}\label{sec:motivation}
Numerous  problems in modern number theory could be solved or at least
better understood, if we have a good knowledge on the  algebraic $K$-theory (or motivic cohomology) of integers of number fields or the cohomology of arithmetic groups (i.e.
 subgroups of finite index of $GL_N(\Z)$). As a short list, we could mention:
\begin{itemize}
\item modular forms and special values of $L$ functions,
\item Iwasawa theory and understanding of the ``cyclotomy'',
\item Galois representations (or automorphic representations).
\end{itemize}
Let us explain first what  is algebraic K-theory: to a commutative ring $R$ we can associate (functorially) an infinite family of abelian groups $K_n(R)$ which encodes a huge amount of information on its arithmetic, geometric and algebraic structures. These groups extend some classical notions and give higher dimensional analogues of some well known results.
For instance if $R$ is a ring of integers or a polynomial ring over a field, we have
\begin{itemize}
\item $K_0(R)$ is the classical Grothendieck group (classifying finitely generated $R-$modules),
\item $K_1(R)$ is the group of invertibles of $R$,
\item $K_2(R)$ classifies the universal extensions of $SL(R)$ and is related to the Brauer group in the case of a field.
\end{itemize}
We can give a general abstract definition of  $K_n$ for $n>0$,
\[
K_n(R)=\pi_n(BGL(R)^+)\, ,
\]
where $BGL(R)^+$ is the Quillen +-construction applied to the classifying space $BGL(R)$ and $\pi_n$ denotes the $n$th homotopy group \cite{Rosenberg1995}.
We also can see the $K$-groups as a way to understand $GL(R)$. For instance, $K_1(R)$ is isomorphic to $GL(R)$ modulo elementary relations. For a more detailed background on K-theory and its applications see \cite{Rosenberg1995}.

{\bfseries Fact}: these groups are hard to compute.%

\subsubsection{The Vandiver conjecture as an illustration}
The Vandiver conjecture plays an important role in the understanding of the ``cyclotomy'' in number theory. Its statement is as follows;\\
\noindent{\bfseries Conjecture of Vandiver}: Let $p$ be an odd prime,
 $\zeta_p=e^{2\pi
i/p}$, $C$ the $p$-Sylow subgroup of the class group of $\Q(\zeta_p)$ and $C^+$ the subgroup fixed by the complex conjugation. The Vandiver conjecture is the statement that $C^+=0$.\\
To illustrate its interplay with K-theory, we have\\
{\bfseries Fact} (Kurihara, 1992)\cite{kurihara}: If $K_{4n}(\Z)=0$ for all   $n>0$, then the conjecture of
 \nobreak{Vandiver}  is true.\\
Some partial results on this conjecture and its connection with the cohomology of $SL_N(\Z)$ are given in \cite{soule-vandiver}.\\
\noindent{\bfseries General problem}: Find explicit methods for computing
(co)homologies of arithmetic groups  and the K-theory of number fields (or their
ring of integers).\\
Our first task will be the computation of the (co)\-ho\-mo\-lo\-gies of linear groups (mainly $GL_N(\Z)$ and $SL_N(\Z)$).
We can show \cite{pfpk_cras,pfpk_mfo,pfpk_egs} that the computation of those groups is the key point for computing $K$-groups with our method. We will begin to recall some facts from topology.
\subsubsection{Topological Excursion: Cellular complexes, or how to simply describe the "combinatorial" structure of a topological space}
The notion of cellular complex is a generalization of a graph to several dimensions. 
We call {\it $n-$cell} a topological space homeomorphic to the open unit ball of $\R^n$
and such that its closure is homeomorphic to the closed unit ball.
A cellular complex (or cell complex or also cellular decomposition or cellular space)  is a family of sets $X^n$ (with  $n\in \N$), such that each $X^n$ is a collection (eventually infinite) of $n-$cells. Usually we work with cell complexes with a finite number of cells.\\
A classical result \cite{Spanier1994} shows that {\itshape any (reasonable) topological space can be approximated by such cell complexes}.
\begin{center}
\begin{tabular}{ccc}
\includegraphics[width=3cm,height=3cm]{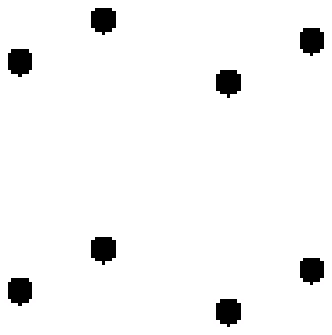} & \qquad\qquad & 
\includegraphics[width=3cm,height=3cm]{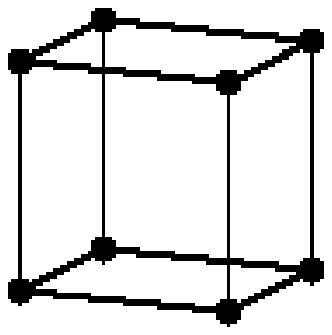} \quad\\[5pt]
\includegraphics[width=3cm,height=3cm]{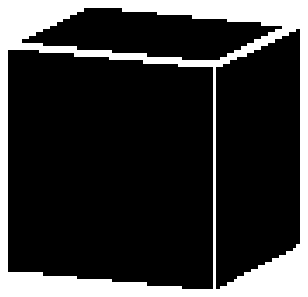}&\qquad\qquad &
\includegraphics[width=3cm,height=3cm]{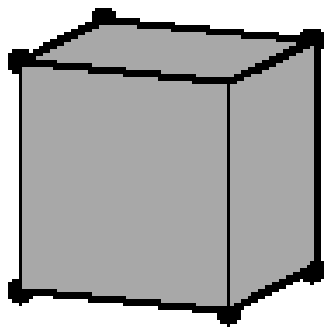}
\end{tabular}
\quad\\
{\bfseries \small A cellular decomposition of the cube}
\end{center}
\subsubsection{Computation of the homology of a cell complex}\label{computation_cell_complex}
To a cell complex, we can associate a family $C_n$ ($n\in \N$) of free $\Z-$modules
and a family $d_n : C_n \to C_{n-1}$ of linear maps.\\
The module $C_n$ is the free module with basis the $n-$cells (modulo a choice of orientation).\\ 
If the complex is finite, then all the modules are of finite rank and we will denote by  \,$\{b_\lambda^n\}_{\lambda \in \Lambda_n}$ a basis of $C_n$, $ \Lambda_n$ being an index set  for the $n-$cells. Then the map $d_n$ is defined by
$$d_n(b_\lambda^n)=\sum_{\mu} [b_\lambda^n : b_\mu^{n-1}] b_\mu^{n-1} \, ,$$
and the integer number $[b_\lambda^n : b_\mu^{n-1}]$ is called {\itshape the incidence number}  of the cell  $e_\mu^{n-1}$ inside the cell $e^n_\lambda$. 
The relation $d_n\circ d_{n+1}=0$ (i.e., formally $d^2=0$) should hold for any $n$.\\ 
If the complex is {\itshape regular} (i.e., always at most one cell of dimension $n+1$ between two cells of dimension $n$), then we can build the incidences inductively starting from the $0-$cells up to the maximal cells using the $d^2=0$ condition and moreover the incidence numbers will be 0, or $\pm1$.\\
The $n$th homology group of the complex is defined as the quotient of $\textrm{Ker}(d_n)$ by $\textrm{Im}(d_{n+1})$. This construction is functorial (in the category of cell complexes).
As a consequence, {\itshape we can determine  the homology groups
effectively by computing the  Smith form of the integral matrices of
the $d_n$ (relatively to the fixed basis).} Notice that the Smith form
gives both the rank of the free part and the explicit description of
the torsion part. In case the computation of the torsion is
unnecessary (or to difficult to achieve), we can tensorize by
$-\otimes_\Z \Q$ and the homology groups become $\Q$-vector spaces
with their dimensions given by the ranks of the matrices of the
differentials. 

In general, the matrices of the differentials can be very large even
for a relatively simple cell complex. However, they are also very
sparse and we may look for some other favorable properties which would enable the
computation despite the size of the problem.

\subsubsection{How to use such settings for the computation of linear groups ?}
If $G$ is a group acting on a cell space $X$ (i.e., $G$ sends
$n-$cells to $n-$cells), then, under some technical assumptions on $X$
and on the action, we can show \cite{B} that roughly  {\itshape computing the
  homology of $G$ (as group homology) is the same as computing the
  homology of the cell space $X/G$}.
Hence, if $X/G$ can be calculated effectively, we can compute explicitly its homology, and from this the homology of $G$ (similarly for the cohomology). Notice that in general the space $X/G$ will not be regular anymore. The main difficulty is to find a cell space $X$ such that $X/G$ will be effective. We will discuss in section \ref{voronoi} how we can construct such cellular space for linear groups.

\subsection{Parallelism motivations}

This first idea to deal with very large sparse matrices is to use them
as blackboxes, i.e. only using the matrix-vector product. This will
let the matrix remain sparse all along the algorithm where Gaussian
elimination for instance would fill it up.
To compute the rank, the fastest black box algorithm is Wiedemann's as
shown in \cite{jgd:2002:villard}.
This algorithm computes a sequence of scalars of the form $u^t A^i v$
($u$ and $v$ are vectors) with $i$ matrix vector products and dot
products.
It has been shown to successfully deal with large sparse
matrix problems from homology, see e.g. \cite{Dumas2001}. 
Nevertheless, when matrices are very large (e.g millions on non-zero entries) computations 
would require months or years. 
This is due to the low practical efficiency of the computation of a sparse matrix-vector product.
For instance, in our case of homology computation, one would need $300$ days of CPU to compute the sequence involving
matrix of $GL_7(\Z)$ with $n=19$ (GL7d19 matrix). To achieve computations of many homologies in a realistic time we then need to parallelize
the computation of the sequence. Then the algorithm candidate is the
block Wiedemann method, which computes a sequence $X^T A^i Y$ where
$X$ and $Y$ are blocks of vectors.
This step can be easily parallelized by distributing vectors of block $Y$ to several 
processors.
We thus have several objectives with regard to the parallelism in this
paper:\\[-.8cm]
\begin{itemize}
\item We want to solve large problems coming from homology computation. 
\item We want to experimentally validate our parallel implementation of the block Wiedemann rank algorithm.
\item We want to show the parallel scaling of block \nobreak{Wiedemann} approaches.
\end{itemize}

\subsection{Summary of the paper}

In section \ref{sec:optim}, the algorithms and optimizations we used
to generate the matrices and compute with them are discussed.
Then, section \ref{sec:result} shows our experimental results on these
large sparse matrices coming from homology.

\section{Algorithms and optimizations}\label{sec:optim}
As seen in section \ref{computation_cell_complex}, we can effectively
compute the homology of a cellular space and of a group which acts
``nicely'' on a cellular space. The main difficulty remains to find
such explicit cellular space. In section \ref{voronoi} we present the
process of matrix generation and optimization. Then in the following
sections we give a description of the algorithms used
for the computation of the rank and the Smith form of those matrices.

\subsection{Matrices generation}\label{voronoi}
In the case of subgroups of  $GL_N(\Z)$, we have an ``obvious'' action on $\Z^N$. We can then capture the topology by regarding $\Z^N$ not as a free $\Z-$module but as a lattice (or equivalently as quadratic forms), and see if this leads to some interesting topological construction. We will describe this approach below and the results that we can get.

\begin{figure*}[t]
{\footnotesize
\[
\setlength{\arraycolsep}{5pt}
\setlength\extrarowheight{2pt}
\begin{array}{|c|c|c|c|c|c|c|c|c|c|c|c|c|c|c|c|c||c|}
\hline
\graycell{\mathbf{n}}  & \graycell{5} & \graycell{6} & \graycell{7} & \graycell{8} & \graycell{9} & \graycell{10} & \graycell{11} & \graycell{12} & \graycell{13} & 
\graycell{14} & \graycell{15} & \graycell{16} & \graycell{17} & \graycell{18} &\graycell{19} & \graycell{20} & \graycell{{\text{{\rm total}}}} \\
\hline
\Sigma_n^*(GL_5(\Z)) & 5& 10& 16& 23& 25& 23& 16& 9& 4& 3 & & & &&& & \mathbf{136}\\
\hline
\Sigma_n(GL_5(\Z)) & 0&0 & 0 & 1 & 7 & 6 & 1& 0 & 2 & 3  & & & & & & & \mathbf{20}\\[1pt]
\hline
\Sigma_n^*(GL_6(\Z) &  3 & 10& 28& 71& 162& 329& 589&  874& 1066& 1039& 775& 425& 181& 57& 18& 7 & \mathbf{5634}\\
\hline
\Sigma_n({ GL_6(\Z)}) & 0 &0 &0 & 0 & 3 & 46 & 163 & 340 & 544 & 636 & 469 & 200 & 49 & 5 & 0 & 0&\mathbf{2455 }\\[1pt]
\hline
\Sigma_n^*(SL_6(\Z)) &  3 & 10 & 28& 71 & 163 & 347 & 691 & 1152 & 1532 & 1551 
& 1134 & 585 & 222 & 62 & 18 & 7& \mathbf{7576}\\
\hline
\Sigma_n({SL_6(\Z)}) & 0 & 3 & 10 & 18 & 43 & 169 & 460 & 815 & 1132 & 1270 
& 970 & 434 & 114 & 27 & 14 & 7& \mathbf{5486}\\[1pt]
\noalign{\hrule height 0.4pt} 

\end{array}
\]
}\caption{Cardinalities of {$\Sigma_n$} and  $\Sigma_n^*$ for $N=5,6$ (empty
slots denote zero)}\label{fig:sig56}
\end{figure*}

\begin{figure*}[t]
{\footnotesize
\[
\setlength{\arraycolsep}{5pt}
\setlength{\extrarowheight}{4pt}
\begin{array}{|c|c|c|c|c|c|c|c|c|c|c|}
\hline
 \graycell{\mathbf{n}}  & & & \graycell{6} & \graycell{7} & \graycell{8} & \graycell{9} & \graycell{10} & \graycell{11} & \graycell{12 }& \graycell{\text{{\rm total}}} \\
\hline
\Sigma_n^*(GL_7(\Z)) &&&{6}& {28} & {115} & {467} & {1882} & {7375} & {26885} & {\mathbf{36758}}\\
\hline
\Sigma_n(GL_7(\Z)) &&&{0}& {0} & {0} & {1} & {60} & {1019} & {8899} & {\mathbf{9979}}\\
\hline\hline
\graycell{\mathbf{n}}  & \graycell{13}& \graycell{14} & \graycell{15} &  \graycell{16} & \graycell{17} & \graycell{18} & \graycell{19} & \graycell{20} & \graycell{21}  & \graycell{\text{{\rm total}}} \\
\hline
\Sigma_n^*(GL_7(\Z)) & {87400}& {244029} & {569568}&  {\mathbf{1089356} }& {\mathbf{1683368} }&  {\mathbf{2075982} }& {\mathbf{2017914} }& {\mathbf{1523376}} & {\mathbf{876385}} &  {\mathbf{10167378}}\\
\hline
\Sigma_n(GL_7(\Z)) & {47271}& {171375} & {460261}&  {\mathbf{955128}} & {\mathbf{1548650}} &  {\mathbf{1955309}} & {\mathbf{1911130}} & {\mathbf{1437547}} & {\mathbf{822922}} &  {\mathbf{9309593}}\\[1pt]
\hline\hline
\graycell{\mathbf{n }} & & & \graycell{22 }& \graycell{23} & \graycell{24 }& \graycell{25 }& \graycell{26} & \graycell{27}& \graycell{\text{{\rm total}}} & {\text{{\rm TOTAL}}}\\
\hline
\Sigma_n^*(GL_7(\Z)) & & &  {374826}& 
 {115411}& {24623}& {3518}& {352}& 33& {\mathbf{518763}} & {\mathbf{10722899}} \\
\hline
\Sigma_n(GL_7(\Z)) & & & {349443}& 
 {105054}& {21074}& {2798}& {305}& 33& {\mathbf{478707}} & {\mathbf{9798279}} \\[1pt]
\noalign{\hrule height 0.4pt} 
\end{array}
\]
}\caption{Cardinalities of {$\Sigma_n$}  and  $\Sigma_n^*$ for $N=7$}\label{fig:sig7}
\end{figure*}

\subsubsection{Vorono{\"\i}'s reduction theory}
Let $N \geqslant 2$ be an integer. We let $C_N$ be the set of positive definite real quadratic forms in $N$ variables. Given $h \in C_N$, let $m(h)$ be the finite set of minimal vectors of $h$, i.e. vectors $v \in {\mathbb Z}^N$, $v \ne 0$, such that $h(v)$ is minimal. A form $h$ is called {\it perfect} when $m(h)$ determines $h$ up to scalar: if $h' \in C_N$ is such that $m(h') = m(h)$, then $h'$ is proportional to $h$. 
\begin{example}
 The form $q(x,y)=x^2+y^2$ has minimum 1 and minimal vectors $\pm (1,0)$ and $\pm (0,1)$. Nevertheless this form is not perfect, because there is an infinite number of  definite positive quadratic forms having these minimal vectors.\\
On the other hand, the form $q(x,y)=x^2+xy+y^2$ has also minimum 1 and has exactly 3 minimal vectors (up to sign), the one above and $\pm (1,-1)$. This form is perfect, the associated lattice is the "honeycomb lattice" (with optimal spheres packing in the plane), it is the only one.
\end{example}
 Denote by $C_N^*$ the set of non negative real quadratic forms on ${\mathbb R}^N$ the kernel of which is spanned by a proper linear subspace of ${\mathbb Q}^N$, by $X_N^*$ the quotient of $C_N^*$ by positive real homotheties, and by $\pi : C_N^* \to X_N^*$ the projection. Let $X_N = \pi (C_N)$ and $\partial X_N^* = X_N^* - X_N$. Let $\Gamma$ be either $GL_N ({\mathbb Z})$ or $SL_N ({\mathbb Z})$. The group $\Gamma$ acts on $C_N^*$ and $X_N^*$ on the right by the formula
$$
h \cdot \gamma = \gamma^t \, h \, \gamma \, , \quad \gamma \in \Gamma \, , \ h \in C_N^* \, ,
$$
where $h$ is viewed as a symmetric matrix and $\gamma^t$ is the transposed of the matrix $\gamma$.
Vorono{\"\i} proved that there are only finitely many perfect forms modulo the action of $\Gamma$ and multiplication by positive real numbers (\cite{Vo}, Th.~p.~110).\\
Given $v \in {\mathbb Z}^N - \{ 0 \}$ we let $\hat v \in C_N^*$ be the form defined by
$$
\hat v (x) = (v \mid x)^2 \, , \ x \in {\mathbb R}^N \, ,
$$
where $(v \mid x)$ is the scalar product of $v$ and $x$. The {\it convex hull} of a finite subset $B \subset {\mathbb Z}^N - \{ 0 \}$ is the subset of $X_N^*$ image by $\pi$ of the elements $\underset{j}{\sum} \, \lambda_j \, \widehat{v_j}$, $v_j \in B$, $\lambda_j \geqslant 0$. For any perfect form $h$, we let $\sigma (h) \subset X_N^*$ be the convex hull of the set $m(h)$ of its minimal vectors. Vorono{\"\i} proved in \cite[\S~8-15]{Vo}, that the cells $\sigma (h)$ and their intersections, as $h$ runs over all perfect forms, define a cell decomposition of $X_N^*$, which is invariant by the action of $\Gamma$. We endow $X_N^*$ with the corresponding $CW$-topology. If $\tau$ is a closed cell in $X_N^*$ and $h$ a perfect form with $\tau \subset \sigma (h)$, we let $m(\tau)$ be the set of vectors $v$ in $m(h)$ such that $\hat v$ lies in $\tau$. Any closed cell $\tau$ is the convex hull of $m(\tau)$ and $m(\tau) \cap m(\tau') = m (\tau \cap \tau')$.

\subsubsection{Vorono{\"\i}'s complex}
Let $d(N) = N(N+1)/2-1$ be the dimension of $X_N^*$ and $n \leqslant d(N)$ a natural integer. We denote by $\Sigma_n^*$ a set of representatives, modulo the action of $\Gamma$, of those cells of dimension $n$ in $X_N^*$ which meet $X_N$, and by $\Sigma_n \subset \Sigma_n^*$ the cells $\sigma$ such that the stabilizer $\Gamma_{\sigma}$ of $\sigma$ in $\Gamma$ preserves its orientation. Let $V_n$ be the free abelian group generated by $\Sigma_n$. We define as follows a map
$$
d_n : V_n \to V_{n-1} \, .
$$

For each closed cell $\sigma$ in $X_N^*$ we fix an orientation of $\sigma$, i.e. an orientation of the real vector space ${\mathbb R} (\sigma)$ of symmetric matrices spanned by the forms $\hat v$, $v \in m(\sigma)$. Let $\sigma \in \Sigma_n$ and let $\tau'$ be a face of $\sigma$. Given a positive basis $B'$ of ${\mathbb R} (\tau')$ we get a basis $B$ of ${\mathbb R} (\sigma)$ by adding after $B'$ a vector $\hat v$, $v \in m(\sigma) - m(\tau')$. We let $\varepsilon (\tau' , \sigma) = \pm 1$ be the sign of the orientation of $B$ in the oriented vector space ${\mathbb R} (\sigma)$ (this sign does not depend on the choice of $v$).

\smallskip

Next, let $\tau \in \Sigma_{n-1}$ be the cell equivalent to $\tau' = \tau \cdot \gamma$. We define $\eta (\tau , \tau') = 1$ (resp. $\eta (\tau , \tau') = -1$) when $\gamma$ is compatible (resp. incompatible) with the chosen orientations of ${\mathbb R} (\tau)$ and ${\mathbb R} (\tau')$.

\smallskip

Finally we define
\begin{equation}
\label{eq1}
d_n (\sigma) = \sum_{\tau \in \Sigma_{n-1}} \sum_{\tau'} \eta (\tau , \tau') \, \varepsilon (\tau' , \sigma) \, \tau \, ,
\end{equation}
where $\tau'$ runs through the set of faces of $\sigma$ which are equivalent to $\tau$.\\
It is shown in \cite{pfpk_egs}, that up to $p$-torsions with $p\le N+1$, the homology of this complex computes the cohomology of $G$.\\
For $N=5,6,7$ we get the following results for $\Sigma_n$.
\begin{theorem} {\rmfamily (Elbaz-Vincent/Gangl/Soul\'e)\cite{pfpk_cras,pfpk_mfo,pfpk_egs}}.\\
The cardinalities of {$\Sigma_n$} and {$\Sigma_n^*$} are shown on figure
\ref{fig:sig56} for $N=5,6$ and  on figure \ref{fig:sig7} for $N=7$.


\end{theorem}
The previous result gives the precise size of the matrices involved in the computation of the homology.

\smallskip

The main challenge was then the computation of the ranks of the matrices of the differential (this gives the free part of the homology) and the computation of the Smith forms (which gives the relevant arithmetical information of the homology), in particular for $N=7$, knowing that such matrices are particularly sparse. We can emphasize the fact that what we want to detect is the ``high torsion'' in the homology (i.e. the prime divisors $>7$ of the Smith invariants).

\smallskip

In the following paragraphs we will discuss the different methods
chosen for the computations and to take up the 
challenge\footnote{All the matrices are available on line in the ``Sparse Integer Matrix
Collection'' (\url{ljk.imag.fr/membres/Jean-Guillaume.Dumas/simc.html})}.

\subsection{Coppersmith Block Wiedemann}

One successful approach to deal with linear algebra computations on
large sparse matrices is to rely on Lanczos/Krylov black-box methods. 
In particular, block versions of Wiedemann method \cite{Wiedemann:1986:SSLE} are well suited for parallel computation.
This technique has been first proposed by Coppersmith in \cite{Coppersmith:1994:SHL} for computation over $GF(2)$ 
 and then analyzed and proved by Kaltofen \cite{Kaltofen:1995:ACB} and Villard \cite{Gilles:further,Gilles:study}.
The idea is for a matrix $A \in \F^{n\times n}$ to compute the minimal generating matrix polynomial of the matrix
sequence $\{XA^iY\}_{i=0}^\infty \in \F^{s \times s}$, where $X,Y$ are
blocks of $s$ vectors (instead of vectors in the original Wiedemann's
algorithm). 
Therefore, an intuitive parallelization of this method is to distribute the 
vectors of the block $X,Y$ to several processors. Thus, the computation of the sequence, which is the major performance 
bottleneck of this method, can be done in parallel and then allow for better performance.

Lots of implementations and practical experimentations has been developed on parallel block Wiedemann.
For instance, in 1996, Kaltofen and Lobo \cite{Kalto-Lobo} have proposed a coarse grain implementation to solve homogeneous linear 
equations over $GF(2)$. 
They have thus been capable to solve a system of $252\,252$ linear equations with about $11.04$ million non-zero entries,
in about $26.5$ hours using 4 processors of an SP-2 multiprocessor. 

Lately, in 2001, Thom\'e in \cite{Tho01} improved Coppersmith's algorithm by introducing matrix half-gcd's computation, and 
its implementation \cite{Thome02b} was able to outperform Kaltofen-Lobo's software. One may remark that introduction of matrix gcd was first 
suggested by Villard in \cite{Gilles:further} who relied on the work of Beckermann and Labahn \cite{BL94} on power Hermite Pad\'e approximation.
Finally, Giorgi, Jeannerod and Villard have generalized in \cite{Giorgi:2003:issac} block Wiedemann algorithms by introducing \sigmabase computation 
and then reducing the complexity to polynomial matrix multiplication.
A sequential implementation of this algorithm is now available in the LinBox library (www.linalg.org).

\subsection{Block symmetry}
In order to reduce the number of dot products, we used
a symmetric projection. In other words, we set $X=Y^T$ in the
$XA^iY$ sequence. Indeed, in this case the probability of success is
reduced but the obtained block is symmetric as soon as $A$ is
symmetric. This is always the case when the preconditioners of
\cite{jgd:2002:villard} are used (they are of the form $A^T A$).
This reduces the dot product part of the computation of the sequence
by a factor of two. 
For instance column $i$ can be deduced from its top $i$ elements and
row $i$. This induces some load balancing issues when one process owns
the computation of one column.
Note also that we use BLAS level-2 for the computation of this dot
products. In other words we perform them by blocks.

\subsection{\sigmabase computations}

In order to efficiently compute \sigmabase we rely on algorithm {\sf PM-Basis} of \cite{Giorgi:2003:issac} which reduces this 
computation to polynomial matrix multiplication. 
One can multiply two polynomial matrices $A,B \in \F^{n \times n}[x]$ of degree $d$ in $O(n^3d+n^2d\log d)$ finite field operations 
if $d$-th primitive roots of unity are available in $\F$. 
Consequently, we decided for our computations to define $\F$ as a prime field with primes of the form $c\times2^k+1$ such
 that $c\times 2^k \equiv 0 \bmod d$. These primes are commonly called FFT primes since they allow the use of FFT evaluation/interpolation 
for polynomials of degree $d$. We refer the reader to \cite{CantorKalto1991, BostanSchost2005} and references therein for 
further informations on fast polynomial matrix arithmetic.

When finite fields not having $d$-th primitive roots of unity are used, polynomial matrix multiplication is still be done efficiently by 
using Chinese Remainder Theorem with few FFT primes. 
Let be $\F$ a prime field of cardinality $p$, then the multiplication of $A,B \in \F^{n \times n}[x]$ of degree $d$ can be efficiently 
done by using CRT with FFT primes $p_i$ satisfying $\prod p_i > d\times n \times p^2$.
This is equivalent to perform the multiplication over the integers and then reduce the result in $\F$.
The overall performance of the multiplication, and then of the \sigmabase, is dependent on the numbers of FFT primes needed.

\subsection{Rank}
Our main interest in the block Wiedemann approach is to compute the rank of large sparse matrices given
 by the Homology group determination problem explained in section \ref{sec:motivation}. 
Hence, we rely on Kaltofen-Saunder's rank algorithm \cite{Kaltofen:1991:SSLS} and its block version \cite{Turner2006} to achieve efficient 
parallel computation.

The Kaltofen-Saunders approach is based on the fact that if $\tilde{A}$ is a good preconditioned matrix of $A$
then its rank is equal to the degree of its minimal polynomial minus
its valuation (or co-degree) \cite{Kaltofen:1991:SSLS}. Thus, by using well chosen preconditioners and Wiedemann algorithm
one can easily compute the rank of a sparse matrix over a finite
field. The block version of this method is presented e.g. in \cite[\S
4]{Turner2006}. We recall now the basic outline of this algorithm.

{\it Block Wiedemann Rank Algorithm} :\\ 
let $A\in\F^{n \times n}$,\\[-.8cm]
\begin{enumerate}
\item[\bf 1] {\bf form} $\tilde{A}$ from $A$ with good preconditioners (e.g. those of \cite{jgd:2002:villard}).
\item[\bf 2] {\bf choose} random block $Y\in\F^{n \times s}$ and {\bf compute } the matrix sequence $S=\{Y^T\tilde{A}^iY\}$ for ${i=0 \hdots 2n/s+O(1)}$.
\item[\bf 3] {\bf compute} the minimal matrix generator $F_{Y}^{\tilde{A}} \in \F^{s \times s}[x]$ of the matrix sequence $S$.
\item[\bf 4] {\bf return} the rank $r$ as $r=\deg(\det(F_{Y}^{\tilde{A}})) - \rm{codeg}   \det(F_{Y}^{\tilde{A}})$.
\end{enumerate}
Note that if the minimal matrix of step 3 is in Popov form
(e.g. computed using the \sigmabase of \cite{Giorgi:2003:issac}),
then the degree of $\det(F_{Y}^{\tilde{A}})$ is
simply the sum of the row degrees of the matrix $F_{Y}^{\tilde{A}}$.
Then the co-degree is zero if the determinant of the constant term of
$F_{Y}^{\tilde{A}}$, seen as a matrix polynomial, is non-zero.
In the latter case the computation of the determinant of the whole
polynomial matrix can be avoided.

When this fails, this determinant is computed by a massively parallel evaluation/interpolation.
It could be interesting, though, to interpolate only the lower coefficients of this polynomial incrementally. 
This was not required for the matrices we considered and we therefore
did not investigate more on these speed improvements. 

Note that to probabilistically 
compute the rank over the integers, it is sufficient to
choose several primes at random and take the largest obtained
value, see e.g. \cite{Dumas2001} for more details. Moreover, one can
choose the primes at random among primes not dividing the determinant
(and thus preserving the rank). In order to ensure this property it it
sufficient to select primes not dividing the valence or last invariant
factor computed by one of the methods of next section.

\subsection{Smith form}

The computation of the Smith form for the matrices of $GL_7(\Z)$ turned out to be a very challenging problem.

\subsubsection{Smith for via the Valence}
Prior experience with sparse homology matrices led us to try the
SmithViaValence algorithm of \cite{Dumas2001}. 
The idea is to compute the minimal valence (the coefficient of the
smallest non zero monomial of the minimal polynomial) of the product $A^TA$ to
determine the primes $p$ which divide the invariant factors of the
Smith form of $A$. 
When the primes have been found, one can compute
the local Smith forms of $A$ at each $p$ separately and return the
resulting Smith form $S$ as the product of the local Smith forms $S_p$
over all $p$. Local Smith form computation can be done by
a repeated Gauss elimination modulo $p^e$ where the exponent $e$ is
adjusted automatically during
the course of the algorithm.

This algorithm works very efficiently for sparse matrices provided that the minimal polynomial of the product $AA^T$ has a
small degree. Unfortunately, the latter condition does not hold in the case of $GL_7(\Z)$
matrices. Moreover, some early experiments with small matrices
showed that much more primes occur in the computed valence than in the
Smith form of the original matrix.

\subsubsection{Saunders and Wan's adaptive algorithm}

Thus, we decided to apply the adaptive algorithm of Saunders and Wan
\cite{Saunders2004} which is a modified version of
Eberly-Giesbrecht-Villard algorithm \cite{Eberly2000}. In
\cite{Eberly2000} the authors proposed a procedure
{\em OneInvariantFactor}($i,A$) (OIF) which computes the $i$th invariant factor of
a $n\times n$ matrix $A$. Then the binary search for distinct invariant factors
allows them to find the Smith form of $A$. OIF reduces the $i$th factor computation to the computation of the last ($n$th)
invariant factor (LIF) of a preconditioned matrix $A+U_iV_i$, where
$U_i,V_i^T$ are random $n\times (n-i)$ matrices. In
\cite{Saunders2004} the method was extended to handle the rectangular
case of $m\times n$ matrix. It is done by computing the last ($i$th) invariant factor
of a preconditioned matrix $L_iAR_i$ where $L_i$ is a $m\times i$ and
$R_i$ is a $i \times n$ matrix. 

The procedure OIF is of Monte Carlo probabilistic type where the
probability of correctness is controlled by repeating the choice of
preconditioners. Assuming the correctness of LIF computation, it
gives a multiple of the $i$th invariant factor. In practice, LIF
is also of randomized Monte Carlo type. The idea is to get a divisor of
the last invariant factor by solving a linear equation $Mx=b$ with random
right-hand side. After several solvings we get the last invariant
factor with large probability. Thus, the overall situation is more complex and we cannot
exclude the possibility that some primes are omitted or unnecessary in
the output of OIF. However, the probability that a prime is omitted or
is unnecessary in this output can be controlled for each prime
separately and is smaller for bigger primes.

Therefore in \cite{Saunders2004} the authors introduce a notion of
smooth and rough parts of the Smith form. 
The idea is to compute the local Smith form
for smaller primes by for example the SmithViaValence or OIF algorithm and to recover
only large
 primes with the invariant factor search of \cite{Eberly2000}. When we
consider large primes, a sufficient probability of correctness can be
obtained by a smaller number of repetitions. 

\subsubsection{More adaptiveness}
As we did not
want to compute the valence, we introduced some minor changes to the
algorithm, which at the end is as follows:
\begin{enumerate}
\item $r = \rank(A)$
\item For primes $1< p < 100$ compute the local Smith form $S_p$ of
$A$;
\item Compute $s_r(A)$ by OneInvariantFactor algorithm;
\item $P$ = all primes $p > 100$ which divide $s_r(A)$;
\item If $P=\emptyset$ return $S=\underset{p}{\Pi} S_p$;
\end{enumerate}

One advantage of this method is that we get the information on
the smooth form of the matrix very quickly. Moreover, the OIF
computation acts as a certification phase which allows us to prove
that no other primes are present with a sufficiently
large probability. This probability is explicit in the following theorem:  
\begin{theorem}\label{thm:oif}
The probability that there exists a prime $p> P$ that divides the
$i$th invariant factor but does not divide the output of OIF which
uses $M$ random preconditioners $L_i, R_i$ and $N$ random vectors $b$,
with $b\in\{0,1,\dots\beta-1\}$ and $\beta>s_i(A)$, in the
LIF procedure is bounded by 
$$
M\sum_{p>P}^{\infty}\left(\frac{2}{p}\right)^N.
$$
\end{theorem}
\begin{proof}
As we take the $\gcd$ of the result with different preconditioners $L_i$
and $R_i$, is suffices that the LIF computation fails in one case to
spoil the computation. We are free to choose a large bound for $\|b\|$
such that $ \|b\|>s_i(A)$ without increasing the complexity of
LIF computation. Then the probability that a prime $p < \|b\|$ is
omitted in LIF is less than or equal to
$\frac{1}{\beta}\lceil\frac{\beta}{p}\rceil< \frac{2}{p}$, see
\cite{Abbott1999}. 
Finally, we bound
the probability that any prime $p>P, p\mid s_i(A)$ is missing by
taking a sum over all primes.
\end{proof}

The choice $M=N=2$ suffice to obtain a small probability $0,015$ of
omitting an important prime, and at the same time to exclude all
primes that are not in the $i$th invariant factor. In our experiments, there was no need to perform the
computation for any additional prime $p > 100$ as all the primes were
excluded by the OIF computation. This is one of the most important
advantages over the valence computation.

 The algorithm \cite{Saunders2004} was stated in the case of dense
 matrices. 
We slightly modified it in order to exploit the sparse structure of the
matrix. In particular, we used the sparse  local Smith
form computation of \cite[Algorithm LRE]{Dumas2001} but stick to the dense Dixon
solver \cite{Dixon1982} as long
as the memory was sufficient. Any other solver,
including the new sparse solver of
\cite{EbGies2006} could potentially be used for larger matrices.

The limits of this method are imposed by the available memory. 
For example, it was possible to use the dense Dixon solver
only for the six smallest matrices. Furthermore, sparse elimination
reached its limits for matrices of size greater than $171375\times
47271$ and $21074\times 105054$ when the filling of the matrices started
to be impossible to handle. For the $460261\times 171375$ matrix GL7d15
and $105054\times 349443$ matrix GL7d23, specialized
space-efficient elimination procedures mod 2, 3 and 5 allowed us to
compute the rank mod 2, 3 and 5 respectively.

\subsubsection{Chain Reductions}

The encountered problems have shown a need for a more elaborated reduction
algorithm. We focused our attention on the algebraic reduction
algorithm for chain complexes of 
\cite{Mrozek1997}. We implemented a simplified version of the
algorithm in the language of matrices using the LinBox library. The
heuristic behind this algorithm is that Gaussian elimination can
propagate from one matrix of a chain complex to the next thanks to the
exactness of the differential map (i.e. $d^2=0$ condition).

\begin{figure}[h]
\begin{center}
\includegraphics[width=6cm]{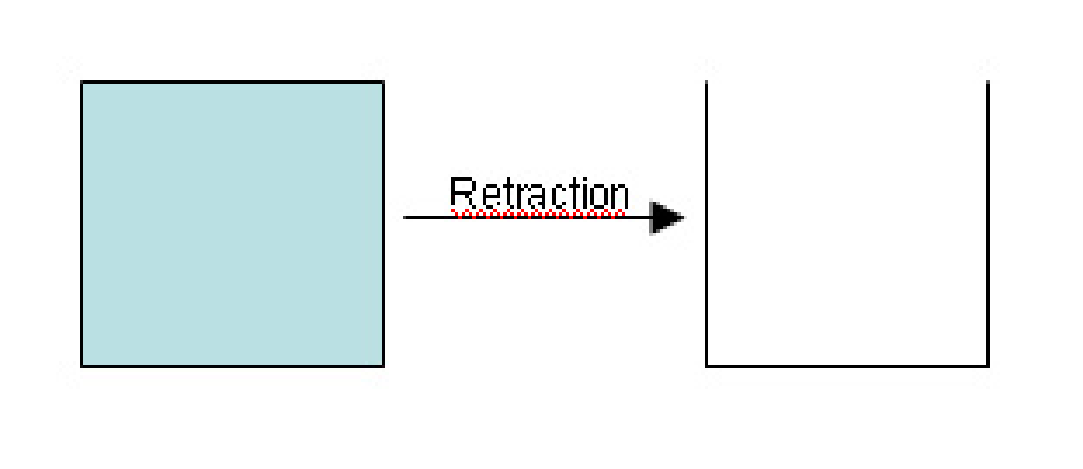}
\end{center}\caption{Retraction for a square}\label{fig:square}
\end{figure}

The motivations come from the geometric properties of homologies. By a
free face we refer to a  $(k-1)$ cell $a$ which is in the differential of
only one $k$-dimensional cell $b$. By removing the pair $(a,b)$ we obtain a retract of the
initial cell complex (viewed as a geometrical object), see
{\bf Figure} \ref{fig:square}. The process can
be repeated. From the homology theory we know that the groups of
homology are the same for the set and its retract. The removal of
pairs leads to a reduction of the basis of the
cell complex. We
refer to \cite[Ch.4]{Mrozek2004} for a full
description of the procedure.  
 
If the differential map is represented by matrices whose rows
represent the cells of dimension $k-1$ and columns - dimension $k$, the removal of
a pair $(a,b)$ can be interpreted as the removal of a row with only
one non-zero entry (1-row) and the column it points to. In the general
case, the algebraic reduction of $(a,b)$ such that $d_k(b)=\lambda a +u$ is possible iff $\lambda$ is
invertible in the ring of computation i.e. $\Q, \Z, \Z_{p^k}$ - depending
on the problem. A modification $\tilde{d}$ of the differential given by the formula
\begin{equation}\label{eq:partial}
\tilde{d}_k v  = d_k(v) - \lambda^{-1}[ v : a]d_k(b).
\end{equation}

Thus, in the basic case of free face removal no modification is needed. In
the case of matrices, the formula describes a step of Gauss
elimination where the reduced row is removed and not permuted. This
proves that the Smith form (or the rank) of the initial and reduced
matrix will be the same, provided we add a number of trivial invariant
factors equal to the number of rows reduced to the reduced Smith form. 

The important characteristic of
this methods is that, thanks to the exactness of the matrix sequence,
we can also remove row $b$ and column $a$ from the neighboring
matrices. In this way, elimination in one matrix can propagate on the
others.

Due to the format of data (large files with matrix entries) we decided
to implement only the simplest case of $1$-rows removal which led to
entries removal but no modifications. We removed the empty
rows/columns at the same time and performed the whole reduction phase
at the moment of reading the files. This led to vast matrix reduction in the
case of $GL_7(\Z)$ matrices from the beginning of the sequence. The
propagation of reductions unfortunately burned out near $GL7d14$
matrix and stopped completely on $GL7d19$. Applying the process for
the transposed sequence did not improve the solution. 
Next step would
be to implement the propagation of Gauss elimination 
steps as in Eq. \eqref{eq:partial}. It would be interesting to examine whether the burn-out can be
connected to the loss of regularity for $X/G$ and/or a huge
rectangularity of the input matrix $GL7d14$.

\section{Experimental results}\label{sec:result}

\subsection{Parallel implementation}

For the sake of simplicity in our experimental validation, we decided to 
develop our parallel implementation using shell tasks distribution on SMP architecture.
Thus, a simple script code is used to distribute all different tasks over all the processors and 
files are used to gather up the computed results.
Our parallel implementation has been done as follow:
\begin{enumerate}
\item The block of vectors $Y$ and the sparse matrix $A$ are broadcasted on every processors.
\item Each process takes one column of the blackbox $A^T.A$ and compute the corresponding column's sequence
   using the first ith columns of $Y$. Each process writes the result in a file labeled with 
the corresponding index of the column sequence.
\item When all previous processes have terminated, the \sigmabase computation is sequentially performed after loading the sequence from the generated files. 
\end{enumerate}

Despite the na\"ive approach used for the parallelization, our
implementation 
authorized us to  perform very large computations as show in next
section. 
However, our experiments show a need for at least 
a more robust parallel computation scheduler.

\subsection{Rank and Smith form}

All our computation have been done on a SGI Altix 3700 gathering 64 Itanium2 processors with $192$Gb memory and 
running SuSE Linux Enterprise System 10. Further informations on this platform are available at 
\url{http://www.math.uwaterloo.ca/mfcf/computing-environments/HPC/pilatus}.

In {\bf Table} \ref{tab:res} we include the information about the
dimensions of the $GL_7(\Z)$ matrices and their sparsity. The matrices
are very sparse which is illustrated by the fact that less than 1\% of
the entries are non-zero except for matrices GL7d10 and GL7d11. This
value drops to less than 0,2\% in the case of the largest matrices.
Also in {\bf Table} \ref{tab:res} we give the results for the rank and the Smith form
computations. We have obtained a full information on 
the rank of $GL_7(\Z)$ matrices. For the computation of the Smith form,
full result has been obtained in the case of matrices
10,11,12,13,25,26. For matrices 14 and 24 only the smooth part of the
Smith form has been computed. For matrices 15 and 23 we have proved
the existence of a non-trivial local Smith form at 2 and 3 and a
triviality of the local Smith form at 5. As the
result for these matrices we give the number of invariant factors
divisible by 2 and 3.

In {\bf Table} \ref{tab:sf-time} we give the times for the Smith form
algorithms used. For cases with * no data are available or
relevant. 

\begin{table*}\footnotesize
\begin{center}
\begin{tabular}{|r|r|r|r|r|r|r|}
\hline 
\graycell{$A$} & \graycell{$\Omega$} & \graycell{$n$} & \graycell{$m$} & \graycell{rank} & \graycell{ker} & \graycell{Smith form} \\
\hline
\hline
GL7d10& 8& 60& 1& 1& 59& 1 \\\hline
GL7d11& 1513&1019& 60& 59& 960& 1 (59) \\\hline
GL7d12& 37519&8899& 1019& 960& 7939& 1 (958), 2 (2)\\\hline
GL7d13& 356232&47271& 8899& 7938& 39333& 1 (7937), 2 (1)\\\hline
GL7d14& 1831183&171375& 47271& 39332& 132043& 1 (39300),2 (29),4 (3)\\\hline
GL7d15& 6080381&460261& 171375& 132043& 28218& 1 (131993), 2$\cdot$? (46),
6$\cdot$? (4) (*)\\\hline
GL7d16& 14488881&955128& 460261& 328218& 626910& \\\hline
GL7d17& 25978098&1548650& 955128& 626910& 921740& \\\hline
GL7d18& 35590540&1955309& 1548650&921740 &1033569 & \\\hline
GL7d19& 37322725&1911130& 1955309&1033568 &877562 & \\\hline
GL7d20& 29893084&1437547& 1911130& 877562& 559985& \\\hline
GL7d21& 18174775&822922& 1437547& 559985& 262937& \\\hline
GL7d22& 8251000&349443& 822922& 262937& 86506& \\\hline
GL7d23& 2695430&105054& 349443& 86505& 18549& 1 (86488), $2\cdot$? (12),
6$\cdot$? (5) (*) \\\hline
GL7d24& 593892& 21074& 105054& 18549& 2525& 1 (18544),2 (4),4 (1)\\\hline
GL7d25& 81671&2798& 21074& 2525& 273& 1 (2507), 2 (18)\\\hline
GL7d26& 7412&305& 2798& 273& 32& 1 (258), 2 (7), 6 (7), 36 (1)
\\\hline
\end{tabular}
\caption{Results of the rank and Smith form computation for $GL_7(\Z)$ matrix
$A$ of dimension $n\times m$ with $\Omega$ non-zero entries. For
(*) the information is incomplete - only divisors of the invariant
factors were determined based on the rank mod 2 and 3 computation.}\label{tab:res}
\end{center}
\end{table*}

\begin{table*}\footnotesize
\begin{center}
\begin{tabular}{|r|r|r|r|r|r|r|r|r|}
\hline
\graycell{$A$}  & \graycell{$\tilde{n}$} & \graycell{$\tilde{m}$} & \graycell{$\tilde{r}$} & \graycell{Red} &
\graycell{RAdaptive} & \graycell{SmoothSF} & \graycell{AdaptiveSF} & \graycell{SFValence}\\
\hline
\hline
GL7d11&39    &8    &52  &0.01s &$<10^{-2}$s&0.09s   &0.26s   & 4.84s\\\hline
GL7d12&289   &58   &909 &0.30s &0.16s  	   &9.75s   &218.68s & 4.04h\\\hline
GL7d13&7938  &740  &7250&3.12s &159.16s	   &0.76h   &* 	     & 2526.65h\\\hline
GL7d14&165450&35741&4279&21.62s& *	   &796h    &*	     & *\\\hline		
GL7d25&2797  &20990&0   &1.74s & *         &17.67s  &4.40h   & 52.13h\\\hline
GL7d26&302   &2748 &0   &0.14s & *         &0.29s   &26.81s  & 274.35s\\\hline
\end{tabular}
\caption{Times for Smith Form computation for $GL_7(\Z)$
matrices. From left to right: the dimensions of the matrix after
reductions, rank approximation by reductions, time of
reading and reducing the matrix, time for the adaptive algorithm for
a reduced matrix; times for the original matrix: smooth form computation, adaptive algorithm; valence
computation in parallel - sequential time equivalent.}\label{tab:sf-time}
\end{center}
\end{table*}

The rank computation for $GL_7(\Z)$ matrices was performed modulo
$65537$. This FFT prime allowed us to use both BLAS routines and
sigma-basis reconstruction using fast polynomial multiplication. In
{\bf Table} \ref{tab:r-opp} we give the timings for different operations used
in the computation i.e. the sparse matrix-vector product, BLAS-based
matrix-vector product and, for the sake of comparison, the time of
scalar dot product equivalent to the BLAS computation. 

\begin{table}[h]\footnotesize
\begin{center}
\begin{tabular}{|r|r|r|r|}
\hline
\graycell{$A$}  & \graycell{$1 A^TAu$ [s]} & \graycell{$1 U^Tv$ [s]} & \graycell{$30 u^Tv$ [s]}\\
\hline 
\hline
GL7d11  &0.0002	&$< 10^{-4}$	&$< 10^{-4}$\\	
GL7d12	&0.0038	&0.0001	&0.0002    \\
GL7d13	&0.0550	&0.0005	&0.0036    \\
GL7d14	&0.2677	&0.0025	&0.0190    \\
GL7d15	&0.9048	&0.0082	&0.0708    \\
GL7d16	&2.7724 &0.0234	&0.2641    \\
GL7d17	&7.8003 &0.0485	&0.5052    \\
GL7d18	&11.9457&0.0759	&0.8710    \\
GL7d19	&13.3591&0.0948	&1.0710    \\
GL7d20  &10.4056&0.0711	&0.8587    \\
GL7d21	&5.7461 &0.0408	&0.4604    \\
GL7d22	&1.9919 &0.0180	&0.2082    \\
GL7d23	&0.4354	&0.0052	&0.0459    \\
GL7d24	&0.0843	&0.0012	&0.0085    \\
GL7d25	&0.0078	&0.0002	&0.0008    \\
GL7d26  &0.0007	& $< 10^{-4}$	&$< 10^{-4}$\\\hline
\end{tabular}\caption{CPU timings (in sec.) for different operations used in large-scale
parallel rank computation. All times in seconds. From left to right:
time of a matrix-vector product, a BLAS multiplication of a vector and a $30
\times \min(n,m)$ matrix $U$ and 30 dot products.}\label{tab:r-opp}
\end{center}
\end{table}

In
{\bf Table} \ref{tab:r-time} we give the estimation of sequential and
parallel cpu time of rank computation and compare it with the real
time of
parallel computation. 
The times are estimated based on the number
of iterations and the times of one step which can be computed from
{\bf Table }\ref{tab:r-opp} (notice, that in the scalar case we use 1 dot product
instead of BLAS). The real time of computation includes the time of
writing and reading the data which was considerable. 
The difference of the real and
estimated running times may also be due to the overload of the
computation cluster. Moreover, long computations suffered from system
crashes and/or shutdowns. Some restoration scripts were used to
recover the data which unfortunately required re-running some part of
the computation. Thus, the real time given in {\bf Table}
\ref{tab:r-time}[Col.5] should be treated as a rough approximation. 
\begin{table*}
\footnotesize
\begin{center}
\begin{tabular}{|r|r|r|r|r|r|r|}
\hline
\graycell{$A$} & \graycell{iter [1]} & \graycell{time app}  & \graycell{iter [$p$]} & \graycell{time}& \graycell{time app} & \graycell{\sigmabase} \\
\hline 
\hline
GL7d11  &120    &0.02s   &6 [30]    &0.01s    &$<10^{-2}$s     &0.58s \\	
GL7d12	&1922   &7.53s   &66 [30]   &0.32s    &0.26s     &12.16s\\
GL7d13	&15878  &880.28s &532 [30]  &51.65s   &29.49s    &249.17s\\
GL7d14	&78666  &5.90h   &2625 [30] &0.56h    &0.20h     &0.45h\\
GL7d15	&264088 &66.98h  &8805 [30] &2.25h    &2.23h     &2.45h\\
GL7d16	&656438 &509.80h &21884 [30]&27.29h   &17.00h    &6.03h\\
GL7d17	&1253822&113.90d &41796 [30]&14d  	&3.80d     &0.57d\\
GL7d18	&1843482&256.50d &46089 [40]&28d  	&6.41d     &1.00d\\
GL7d19	&2067138&321.89d &41345 [50]&35d  	&6.44d     &1.56d\\
GL7d20  &1755126&212.82d &36568 [48]&10d  	&4.43d     &1.41d\\
GL7d21	&1119972&75.01d  &37335 [30]&5d   	&2.50d     &0.55d\\
GL7d22	&525876 &293.60h &17532 [30]&16.47h   &9.79h     &5.85h\\
GL7d23	&173012 &21.18h  &5769 [30] &1.17h     &0.71h	    &1.09h\\
GL7d24	&37100  &3172.79s&1239 [30] &188.78s  &105.96s   &666.83s\\
GL7d25	&5052   &40.21s  &171 [30]  &1.56s    &1.36s     &41.47s\\
GL7d26  &548    &0.40s   &21 [30]   &0.03s    &0.02s     &2.03s\\\hline
\end{tabular}\caption{A summary of large-scale parallel rank
computation. From left to right: number of iteration in the scalar
case, time estimation in this case, number of iterations on $p$
processors computed as $2+2\cdot r/p$, average (real) time of sequence
computation, estimated time on $p$ processors,
the time of the $\sigma$ basis computation for a sequence of length
iter [$p$] of $p \times p$ matrices.}\label{tab:r-time}
\end{center}
\end{table*}

\section{Conclusion}
Using the previous methods and computations, we get the following new result for the rational cohomology of $GL_7(\Z)$.
\begin{theorem} {\rmfamily (Elbaz-Vincent/Gangl/Soul\'e)\cite{pfpk_egs}}
 We have
$$
H^m ({GL}_7 (\Z) , \Q) = \begin{cases} \Q &{\rm if }\quad m=0, 5, 11, 14, 15,\\
0 & {\rm otherwise}.
\end{cases}
$$
\end{theorem}

Clearly the simple parallelization we used together with the highly
optimized routines were the key to enable these computations.

To go further and solve even larger problems, it is mandatory to improve
the parallelism. On SMP we can split the matrices into blocks and
perform the matrix-vector products with different
threads. This, and the unbalanced load we had when we choose to assign
one vector to one process, advocates for the use of more advanced
scheduling. We are experimenting KAAPI\footnote{Kernel for Adaptive,
  Asynchronous Parallel and Interactive programming,
  \url{kaapi.gforge.inria.fr}} but were not ready for the
computation of $GL_7$. 

Other improvements are of algorithmic type. They include the use of
the sparse projections of \cite{EbGies2006} for the matrix
sequence. But then we loose the symmetry of the projections and
therefore must pay a factor of two for the number of iterations. 
We could also use an early termination strategy to stop the iteration earlier,
but up to now this require to loose the fast algorithm for the sigma
bases. 
Then if a good structure for the sparsity of the matrices could be
found, e.g. an adapted reordering technique, this would
enable an efficient clustering and therefore faster and more scalable
matrix-vector products.

In order to have the relevant part of the torsion of the integral
cohomology of ${GL}_7 (\Z)$, we would need the complete description of
the Smith forms of all the matrices described above. Our experiments
have shown that this can be an enormously difficult task. The
computation remains to be done but would have applications in number
theory.
\section{Acknowledgments}
 We are grateful to Arne Storjohann and to the Computer Science Computing Facilities of the University of Waterloo
 for letting us fill up their SMP machine to perform our parallel computations.

\bibliographystyle{abbrv}
\bibliography{parank.bib}

\end{document}